\documentclass[a4paper,UKenglish,cleveref,autoref,thm-restate]{oasics-v2021}

\title{Non-deterministic updates of Boolean networks} 

\author{Loïc Paulevé}{Univ. Bordeaux, Bordeaux INP, CNRS, LaBRI, UMR5800, F-33400 Talence, France
\and
\url{https://loicpauleve.name}}{loic.pauleve@labri.fr}{https://orcid.org/0000-0002-7219-2027}{French
Agence Nationale pour la  Recherche  (ANR): ANR-FNR project ``AlgoReCell'' ANR-16-CE12-0034 and ANR
project ``BNeDiction'' ANR-20-CE45-0001.}
\author{Sylvain Sené}{Université publique, Marseille, France
\and \url{http://pageperso.lis-lab.fr/~sylvain.sene/}}{sylvain.sene@lis-lab.fr}{}%
{his salary as a French state agent affiliated to Université Aix-Marseille, Université Toulon, CNRS, LIS, UMR 7020, Marseille, France, and the ANR project ``FANs'' ANR-18-CE40-0002.}

\authorrunning{L. Paulevé and S. Sené}
\Copyright{Loïc Paulevé and Sylvain Sené}



\category{} 

\relatedversion{} 

\nolinenumbers 

\hideOASIcs 

\EventEditors{Alonso Castillo-Ramirez, Pierre Guillon, and K\'{e}vin Perrot}
\EventNoEds{3}
\EventLongTitle{27th IFIP WG 1.5 International Workshop on Cellular Automata and Discrete Complex Systems (AUTOMATA 2021)}
\EventShortTitle{AUTOMATA 2021}
\EventAcronym{AUTOMATA}
\EventYear{2021}
\EventDate{July 12--14, 2021}
\EventLocation{Aix-Marseille University, France}
\EventLogo{}
\SeriesVolume{90}
\ArticleNo{10}

\usepackage{stmaryrd}
\usepackage{tikz}

\newcommand{\trans}{\xrightarrow}

\newcommand{\one}{\textsf{1}}
\newcommand{\zero}{\textsf{0}}
\newcommand{\B}{\mathbb B}
\newcommand{\N}{\mathbb N}

\newcommand{\flip}[2]{{#1}^{\overline #2}}
\newcommand{\range}[1]{\llbracket {#1}\rrbracket}

\DeclareMathOperator{\D}{\mathscr{D}}
\DeclareMathOperator{\e}{\mathsf{e}} 
\DeclareMathOperator{\p}{\mathsf{p}} 
\DeclareMathOperator{\bs}{\mathsf{bs}} 
\DeclareMathOperator{\per}{\mathsf{per}} 
\DeclareMathOperator{\ga}{\mathsf{a}} 
\DeclareMathOperator{\fa}{\mathsf{fa}} 

\DeclareMathOperator{\M}{\mathsf{M}}
\DeclareMathOperator{\Mb}{\mathsf{\overline M}}
\DeclareMathOperator{\MP}{\mathsf{MP}}

\DeclareMathOperator{\interval}{Int}
\newcommand{\PhiD}[2]{\Phi_{#2/#1}}
\newcommand{\PhiI}[1]{\Phi_{\interval,#1}}

\renewcommand{\land}{\text{ and }}

\begin{document}

\maketitle

\begin{abstract}
Boolean networks are discrete dynamical systems where each automaton has its own Boolean function for computing its state according to the configuration of the 
network.
The updating mode then determines how the configuration of the network evolves over time.
Many of updating modes from the literature, including synchronous and asynchronous modes, can be defined as the composition of elementary deterministic configuration 
updates, i.e., by functions mapping configurations of the network.
Nevertheless, alternative dynamics have been introduced using ad-hoc auxiliary objects, such as that resulting from binary projections of Memory Boolean networks, or that resulting from additional pseudo-states for Most Permissive Boolean networks.
One may wonder whether these latter dynamics can still be classified as updating modes of finite Boolean networks, or belong to a different class of dynamical systems.
In this paper, we study the extension of updating modes to the composition of non-deterministic updates, i.e., mapping sets of finite configurations.
We show that the above dynamics can be expressed in this framework, enabling a better understanding of them as updating modes of Boolean networks.
More generally, we argue that non-deterministic updates pave the way to a unifying framework for expressing complex updating modes, some of them enabling transitions that cannot be computed with elementary and non-elementary deterministic updates.
\end{abstract}


\section{Introduction}

Boolean networks (BNs) are formal dynamical systems composed of automata, each of them having a Boolean state.
A major difference between BNs and cellular automata (CAs) is that each automaton of a BN follows its own rules for computing its next state depending on the states of the other automata in the network.
Consequently, whereas influences between cells in a CA are structured homogeneously according to a
cellular space, those between automata in a BN are structured according to any directed graph.
In this paper, only finite BNs are considered, as it is generally the case in the literature,
notably because BNs are mostly viewed as both a real-world computational model and a real-world
modeling framework. 

The study of BNs led to fundamental results linking the network architecture (structure of influences between automata) to the existence of fixed points and to 
the number of limit cycles they can exhibit~\cite{J-Aracena2008,J-Demongeot2012,C-Bridoux2021}.
Notably, it is well known that such limit behaviors may depend on the way automata update their state over time~\cite{J-Aracena2009,T-Elena2009,J-Aracena2011,T-Noual2012}.
This emphasizes the importance of what is classically called the updating modes in the analyses of BNs.

BNs are widely employed to model natural systems, with prominent applications in biology.
These applications inspired the definition of various updating modes aiming at reflecting
constraints related to the quantitative nature of the abstracted system, such as reaction duration
and influence thresholds.
There is actually no consensus about one updating mode that would be the most likely, the most representative of the biological reality. 
As a consequence, the choice of this or that updating mode strongly depends on the problematics, on the nature of the questions addressed.
Thus, it remains essential to analyze the impact of a wide range of updating modes with distinct features. 

In this paper, we address the formalization of updating modes in the framework of BNs.
From a very general perspective, given a BN and one of its configurations, an updating mode specify
how to compute the possible next configurations (plural implying non-deterministic systems).

A large majority of updating modes introduced so far can be expressed using deterministic functions
mapping the configurations of the network.
This leads to \emph{elementary} transitions, as it is the case with synchronous (or parallel) and
asynchronous~\cite{J-Noual2018} updating modes, which may result in non-deterministic dynamics.
These functions may also be composed, as in block-sequential~\cite{B-Robert1986} and
block-parallel~\cite{J-Demongeot2020} updating modes, generating non-elementary transitions.

These compositions of deterministic \emph{updates}, however, do not cover all the updating modes introduced in the literature.
Indeed, updating modes may also make use of parameters that cannot \emph{a priori} and intuitively be directly captured by these deterministic updates.
These parameters can represent kinds of delays or threshold effects of state changes.
In this paper, we focus on 3 examples of BN dynamics which have been recently introduced and defined using ad-hoc
formalizations:
\begin{itemize}
    \item \emph{Memory Boolean networks} (MBNs)~\cite{J-Goles2020,C-Goles2019} take into account
        some kind of delay for the decrease of automata.  They have been introduced by the means of
        a deterministic dynamical system with non-binary configurations, whose updates are computed
        deterministically from the BN and a memory vector, specifying the delay for each automaton.

    \item \emph{Interval Boolean networks} (IBNs)~\cite{beyond-general} account for a duration for
        updating an automaton. The other automata can be updated until the former
        automaton eventually change of state.
        They have been defined by an encoding as the fully-asynchronous updating of a BN of
        dimension $2n$.
        The dynamics of the original BN are then recovered by projection.

    \item \emph{Most Permissive Boolean networks} (MPBNs)~\cite{MPBNs} bring a formal abstraction of
        trajectories of quantitative models which are compatible with the BN formalism: from an initial
        configuration, if there is no trajectory where a given automaton is 1 (or 0), then, no
        quantitative refinement of the model can increase (or decrease) the value of this automaton.
        MPBNs have been defined by introducing additional states for automata to account for their
        state change (increasing and decreasing). An automaton in one of these states can be read
        non-deterministically as 0 or 1.
\end{itemize}
Overall, the definition of these BN dynamics involve either non-Boolean configurations, projections
of higher-dimension BN, or both.
Importantly, they suggest that deterministic updates are not expressive enough to capture
specific dynamics.
This is striking with IBNs and MPBNs which can generate transitions that are neither elementary nor
non-elementary transitions, and thus predict trajectories that are impossible with the asynchronous
updating mode.

We show that these dynamics can all be expressed using Boolean configurations in a simple generic framework, which extends the deterministic updates to \emph{non-deterministic} updates:
functions mapping sets of configurations.
In the case of MBNs, the obtained definition from the binary projections of their deterministic
discrete dynamics actually help to understand the generated dynamics:
the transitions match with a particular subset of elementary transitions, suggesting a simpler
parameterization.
In the case of IBNs and MPBNs, the transitions extend the elementary and non-elementary transitions
by considering some delay for the state changes, and having different interpretation of how to ``read''
an automaton in the course of state change.
The obtained definitions suggest many variants for generating sub-dynamics, similarly to the asynchronous mode which generates all elementary transitions.

Thus, non-deterministic updates offer a unified yet simple framework for defining and understanding
BN updating modes with more expressivity than usual deterministic updates.
However, should any set update be considered as a BN updating mode?
We propose an argumentation for a \emph{reasonable} updating mode in the last section, where we suggest that the state change should always be justified by the application of a local function.
This suggests that the MP updating mode generates the largest set of transitions that fulfill this
criterion.

\subparagraph*{Notations}
The Boolean domain $\{0,1\}$ is denoted by $\B$; the set $\{1,\cdots,n\}$ is denoted by $\range n$.
Given a finite domain $A$ with a partial order $\preceq$, and a function $h$ mapping elements of $A$ to $A$,
for any $k\in\mathbb N_{>0}$, we write $h^k$ for $h$ iterated $k$ times.
Whenever for any $a\in A$, $a \preceq h(a)$, we write $h^\omega$ for the iteration of $h$ until
reaching a fixed point
(in this paper, $A$ is often a power set with $\preceq$ being the subset relation).

\section{Boolean networks and dynamics}

A \emph{Boolean network} (BN) of dimension $n$ is specified by a function $f: \B^n \to \B^n$ mapping Boolean vectors of dimension $n$.
The components $\range n$ of the BN are called \emph{automata}.
For each automaton $i \in \range n$, $f_i:\B^n\to \B$ is the $i$-th component of this function, that we call the \emph{local function} of automaton $i$.
The $2^n$ Boolean vectors of $\B^n$ are called the \emph{configurations} of the BN.
In a configuration $x\in\B^n$, $x_i$ is the \emph{state} of automaton $i$.

\subparagraph*{Updating modes}
Given a BN $f$ of dimension $n$ and one of its configurations $x\in\B^n$, an \emph{updating mode} $\mu$
characterizes the possible evolutions of $x$ with respect to $f(x)$.
The dynamical system $(f, \mu)$  defines a binary \emph{transition relation} between configurations of $\B^n$ denoted by $\trans{}_{(f,\mu)}\ \subseteq\ \B^n\times\B^n$.
This dynamical system can be represented by a directed graph $\D_{(f,\mu)} = (\B^n,
\trans{}_{(f,\mu)})$.
This graph is usually called the \emph{transition graph} of $(f, \mu)$.
The reflexive and transitive closure of relation $\trans{}_{(f,\mu)}$, denoted by
$\trans{}_{(f,\mu)}^*$ can be defined as follows: given two configurations $x, y \in \B^n$, $x \trans{}_{(f,\mu)}^* y$ if and only if $x=y$ or there exists a path from $x$ to $y$ in $\D_{(f,\mu)}$.

A \emph{deterministic} updating mode ensures that, for any BN $f$ of dimension $n$, each configuration has at most one outgoing transition ($\forall x,y,z\in\B^n$, $x\trans{}_{(f,\mu)} y$ and $x\trans{}_{(f,\mu)} z$ only if $y=z$).
Otherwise, the updating mode is qualified as \emph{non-deterministic}.

\smallskip

In the following, we consider the BN $f$ to be fixed, and thus, for the sake of simplicity, we omit
the subscript $f$: the transition relation is denoted by $\trans{}_\mu$ and the transition
graph by $\D_\mu$.

\subparagraph*{Dynamical properties}
A configuration $x \in \B^n$ is \emph{transient} if there exists a configuration $y$ such that $x \trans{}_{\mu}^* y$ and $y \not\trans{}_{\mu}^* x$.
Configurations that are not transient are called \emph{limit configurations}.
Because $n$ is finite, these configurations induce the terminal strongly connected components of $\D_\mu$, called the \emph{limit sets} of $(f,\mu)$.
If there exists at least one path from a transient configuration to a limit set, this limit set is called an \emph{attractor} of $(f,\mu)$~\cite{C-Cosnard1985,J-Milnor1985}.
The \emph{basin of attraction} of an attractor $\mathcal{A}$ of $(f,\mu)$, denoted by $\mathcal{B}(\mathcal{A})$, is the sub-graph of $\D_{\mu}$ induced by the set of transient configurations $x$ such that, for any limit configuration $y$ belonging to $\mathcal{A}$, $x \trans{}_{\mu}^* y$.
A limit set of cardinal $1$, \emph{i.e.} composed of a unique limit configuration $x$ is called a \emph{fixed point} of $(f,\mu)$.
A limit set of cardinal greater than $1$ is called a \emph{limit cycle} of $(f,\mu)$.

\section{Updating modes with deterministic updates}

\paragraph*{Elementary transitions}

Let us consider a BN $f$ of dimension $n$ and one of its configurations $x\in\B^n$.
Whenever $x$ and $f(x)$ differ by more than one component, one may define several ways to update $x$:
either by replacing it with $f(x)$, i.e., applying simultaneously the local functions on every automata, or by modifying the state of only a subset of automata.
For each set of automata to update, we obtain a deterministic function mapping configurations, that we refer to as an \emph{elementary} deterministic update:
\begin{definition}
Given a BN $f$ of dimension $n$ and a set of automata $W\subseteq \range n$, $\phi_W : \B^n\to \B^n$ is an \emph{elementary deterministic update} with
\begin{equation*}
    \forall x \in \B^n, \forall i \in \range n,\ \phi_W(x)_i = \begin{cases}
        f_i(x) & \text{if } i \in W\text{,}\\
        x_i & \text{otherwise.}
    \end{cases}
\end{equation*}
\end{definition}
Whenever referring to singleton sets $\{i\}$ with $i\in\range n$, we write $\phi_i$ instead of $\phi_{\{i\}}$.
Notice that $\phi_{\range n} = f$.

\begin{example}\label{ex:Gf}
	Let us consider the BN $f$ of dimension $n = 3$ with 
    $f(x)=\left(\begin{array}{l}
    f_1(x) = \neg x_3 \\ f_2(x) = \neg x_1 \wedge x_3 \\ f_3(x) = \neg x_1
\end{array}\right)$.

    \autoref{tab:example} shows four distinct updatings on its configurations.
The first updating is ineffective and consists in changing nothing.
The second updating changes the state of automaton $1$ by application of $\phi_1$, the third one changes the states of both automata $2$ and $3$ by application of $\phi_{\{2,3\}}$, and the fourth one changes the state of every automaton by application of $\phi_{\range 3}$.

\begin{table}[t!]
	\caption{Configurations, local functions ($\left(f_i\right)_{i \in \range 3}$) and four updating functions ($\phi_\emptyset$, $\phi_1$, $\phi_{\{2,3\}}$, and $\phi_{\range 3}$) of Boolean network $f$ presented in Example~\ref{ex:Gf}.\label{tab:updates}}{%
	\centerline{\scalebox{.85}{\begin{tabular}{@{}c|ccc|cccc@{}}
		\hline
		$x = (x_1, x_2, x_3)$ &
		$f_1(x)$ & 
		$f_2(x)$ & 
		$f_3(x)$ & 
		$\phi_\emptyset(x)$ & 
		$\phi_1(x)$ & 
		$\phi_{\{2,3\}}(x)$ & 
		$\phi_{\range 3}(x) \equiv f(x)$\\
		\hline
		 $(\zero,\zero,\zero)$ &
		$\one$ &
		$\zero$ &
		$\one$ &
		$(\zero,\zero,\zero)$ &
		$(\one,\zero,\zero)$ &
		$(\zero,\zero,\one)$ &
		$(\one,\zero,\one)$ \\
		 $(\zero,\zero,\one)$ &
		$\zero$ &
		$\one$ &
		$\one$ &
		$(\zero,\zero,\one)$ &
		$(\zero,\zero,\one)$ &
		$(\zero,\one,\one)$ &
		$(\zero,\one,\one)$ \\
		 $(\zero,\one,\zero)$ &
		$\one$ &
		$\zero$ &
		$\one$ &
		$(\zero,\one,\zero)$ &
		$(\one,\one,\zero)$ &
		$(\zero,\zero,\one)$ &
		$(\one,\zero,\one)$ \\
		 $(\zero,\one,\one)$ &
		$\zero$ &
		$\one$ &
		$\one$ &
		$(\zero,\one,\one)$ &
		$(\zero,\one,\one)$ &
		$(\zero,\one,\one)$ &
		$(\zero,\one,\one)$ \\
		 $(\one,\zero,\zero)$ &
		$\one$ &
		$\zero$ &
		$\zero$ &
		$(\one,\zero,\zero)$ &
		$(\one,\zero,\zero)$ &
		$(\one,\zero,\zero)$ &
		$(\one,\zero,\zero)$ \\
		 $(\one,\zero,\one)$ &
		$\zero$ &
		$\zero$ &
		$\zero$ &
		$(\one,\zero,\one)$ &
		$(\zero,\zero,\one)$ &
		$(\one,\zero,\zero)$ &
		$(\zero,\zero,\zero)$ \\
		 $(\one,\one,\zero)$ &
		$\one$ &
		$\zero$ &
		$\zero$ &
		$(\one,\one,\zero)$ &
		$(\one,\one,\zero)$ &
		$(\one,\zero,\zero)$ &
		$(\one,\zero,\zero)$ \\
		 $(\one,\one,\one)$ &
		$\zero$ &
		$\zero$ &
		$\zero$ &
		$(\one,\one,\one)$ &
		$(\zero,\one,\one)$ &
		$(\one,\zero,\zero)$ &
		$(\zero,\zero,\zero)$ \\
		\hline
	\end{tabular}}}}{}
    \label{tab:example}
\end{table}
\end{example}

We can then define the notion of \emph{elementary transitions} of a BN, that are the transitions obtained by applying any elementary update on a non-empty subset of automata.
\begin{definition}
    Given a BN $f$, its \emph{elementary} transitions $\trans{}_{\e}\ \subseteq\,\B^n\times\B^n$ are such that, for all configurations $x,y\in\B^n$, $x\trans{}_{\e} y$ if and only if there exists a non-empty subset of automata $W\subseteq\range n$ with $y=\phi_W(x)$.
\end{definition}

Let us now define some classical deterministic and non-deterministic updating modes from these elementary updates.

\paragraph*{Examples of deterministic updating modes}

The most direct updating mode is the application of $f$ to the configuration $x$, resulting in the configuration $f(x)$, or, equivalently, $\phi_{\range n}(x)$:
\begin{definition}
    The \emph{synchronous} (or \emph{parallel}) updating mode of a BN $f$ of dimension $n$ generates the transition relation $\to_{\p}\ \subseteq\,\B^n\times\B^n$ such that, for all configurations $x,y\in\B^n$, $x\to_{\p} y$ if and only if $y=\phi_{\range n}(x)$.
\end{definition}

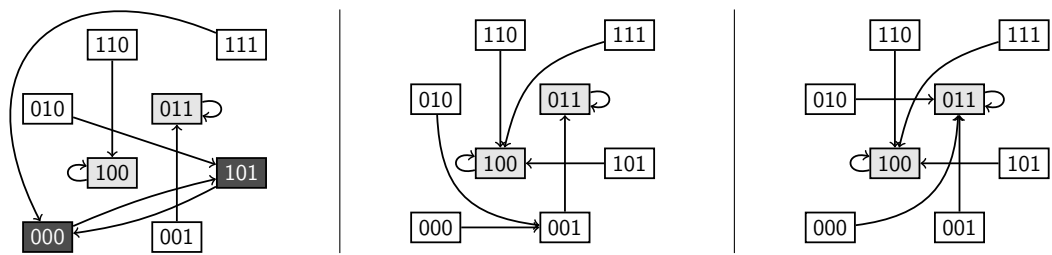
\begin{figure}[t!]
	\centerline{
  \begin{minipage}{.25\textwidth}
	  \centerline{\scalebox{.85}{\begin{tikzpicture}[>=to,auto]
\path[use as bounding box] (-0.6,-0.25) rectangle (3.4,3.55);
			\tikzstyle{type} = []
			\tikzstyle{conf} = [rectangle, draw, thick]
			\tikzstyle{pf} = [rectangle, draw, thick, fill=black!10]
			\tikzstyle{lc} = [rectangle, draw, thick, fill=black!70]
			\node[lc](n000) at (0,0) {\textcolor{white}{$\zero\zero\zero$}};
			\node[conf](n001) at (2,0) {$\zero\zero\one$};
			\node[conf](n010) at (0,2) {$\zero\one\zero$};
			\node[pf](n011) at (2,2) {$\zero\one\one$};
			\node[pf](n100) at (1,1) {$\one\zero\zero$};
			\node[lc](n101) at (3,1) {\textcolor{white}{$\one\zero\one$}};
			\node[conf](n110) at (1,3) {$\one\one\zero$};
			\node[conf](n111) at (3,3) {$\one\one\one$};
        \draw [thick, ->] (n000) edge[bend left=5] (n101);
		\draw [thick, ->, bend left=10] (n101) edge (n000);
			\draw[thick, ->] (n100) edge[loop left,distance=4mm] (n100);
			\draw[thick, ->] (n110) edge (n100);
			\draw[thick, ->] (n010) edge (n101);
			\draw[thick, ->] (n011) edge[loop right,distance=4mm] (n011);
		\draw [thick, ->, bend right=70, looseness=1.75] (n111) edge (n000);
			\draw[thick, ->] (n001) edge (n011);
		\end{tikzpicture}}}
  \end{minipage}
  \hfill\vrule\hfill
  \begin{minipage}{.25\textwidth}
	  \centerline{\scalebox{.85}{\begin{tikzpicture}[>=to,auto]
			\tikzstyle{type} = []
			\tikzstyle{conf} = [rectangle, draw, thick]
			\tikzstyle{pf} = [rectangle, draw, thick, fill=black!10]
			\tikzstyle{lc} = [rectangle, draw, thick, fill=black!70]
			\node[conf](n000) at (0,0) {$\zero\zero\zero$};
			\node[conf](n001) at (2,0) {$\zero\zero\one$};
			\node[conf](n010) at (0,2) {$\zero\one\zero$};
			\node[pf](n011) at (2,2) {$\zero\one\one$};
			\node[pf](n100) at (1,1) {$\one\zero\zero$};
			\node[conf](n101) at (3,1) {$\one\zero\one$};
			\node[conf](n110) at (1,3) {$\one\one\zero$};
			\node[conf](n111) at (3,3) {$\one\one\one$};

			\draw[thick, ->] (n000) edge (n001);
			\draw[thick, ->] (n001) edge (n011);
			\draw[thick, ->] (n100) edge[loop left,distance=4mm] (n100);
			\draw[thick, ->] (n110) edge (n100);
		\draw [thick, ->, bend right=40, looseness=1.25] (n010) edge (n001);
			\draw[thick, ->] (n011) edge[loop right,distance=4mm] (n011);
		\draw [thick, ->, bend right, looseness=1.2] (n111) edge (n100);
			\draw[thick, ->] (n101) edge (n100);
		\end{tikzpicture}}}
  \end{minipage}
  \hfill\vrule\hfill
  \begin{minipage}{.25\textwidth}
	  \centerline{\scalebox{.85}{\begin{tikzpicture}[>=to,auto]
			\tikzstyle{type} = []
			\tikzstyle{conf} = [rectangle, draw, thick]
			\tikzstyle{pf} = [rectangle, draw, thick, fill=black!10]
			\tikzstyle{lc} = [rectangle, draw, thick, fill=black!70]
			\node[conf](n000) at (0,0) {$\zero\zero\zero$};
			\node[conf](n001) at (2,0) {$\zero\zero\one$};
			\node[conf](n010) at (0,2) {$\zero\one\zero$};
			\node[pf](n011) at (2,2) {$\zero\one\one$};
			\node[pf](n100) at (1,1) {$\one\zero\zero$};
			\node[conf](n101) at (3,1) {$\one\zero\one$};
			\node[conf](n110) at (1,3) {$\one\one\zero$};
			\node[conf](n111) at (3,3) {$\one\one\one$};

		\draw [thick, ->, bend right=40, looseness=1.2] (n000) edge (n011);
			\draw[thick, ->] (n011) edge[loop right,distance=4mm] (n011);
			\draw[thick, ->] (n100) edge[loop left,distance=4mm] (n100);
			\draw[thick, ->] (n110) edge (n100);
			\draw[thick, ->] (n010) edge (n011);
		\draw [thick, ->, bend right, looseness=1.2] (n111) edge (n100);
			\draw[thick, ->] (n101) edge (n100);
			\draw[thick, ->] (n001) edge (n011);
		\end{tikzpicture}}}
  \end{minipage}}
	\caption{Distinct possible block-sequential dynamics of BN $f$ defined in Example~\ref{ex:Gf}: 
		(left panel) its parallel dynamics associated with ordered partition $(\range 3)$; 
		(central panel) the block-sequential dynamics associated with $(\{2,3\},\{1\})$;
		(right panel) the sequential dynamics associated with $(\{3\},\{1\},\{2\})$.}
	\label{fig:bs}
\end{figure}

\emph{Sequential} updating modes are parameterized by a permutation of $\range n$, fixing an ordering of elementary updates of single automata~\cite{J-Fogelman1983,J-Goles2008,J-Demongeot2010}.
They can be generalized to \emph{block-sequential} updating modes~\cite{B-Robert1986,J-Aracena2009,C-Goles2010}, parameterized by a permutation of a partition of $\range n$:
\begin{definition}
    Given a BN $f$ of dimension $n$ and $\bs=(W_1,\cdots,W_p)$ an ordered partition of $\range n$, the \emph{block-sequential} updating mode generates the transition relation $\to_{\bs}\ \subseteq\,\B^n\times\B^n$ such that, for all configurations $x,y\in\B^n$, $x\to_{\bs} y$ if and only if $y=\phi_{W_p}\circ\cdots \phi_{W_1}(x)$.
\end{definition}

Remark that the transitions of sequential and block-sequential modes may not be elementary. 
However, they always correspond to a path of elementary transitions: 
$x\to_{\bs} y$ only if $x\to_{\e}^* y$.

Going further in generalization, one may consider deterministic updating modes as infinite sequences of sets of automata, so that automata of a same subset execute their local function in parallel while the subsets are iterated sequentially.
Remark that any of these possible deterministic updating modes will generate transitions corresponding to specific paths of elementary transitions.


\paragraph*{Examples of non-deterministic updating modes}

It is important to notice that deterministic updates can lead to non-deterministic dynamics by allowing different updates on a same configuration.
The most obvious example is the \emph{asynchronous} mode\footnote{The \emph{asynchronous} mode is often referred to as \emph{general asynchronous} in the systems biology modeling community.} consisting of all the elementary transitions.
\begin{definition}
    The \emph{asynchronous} updating mode of a BN $f$ generates the transition relation $\to_{\ga}\ \subseteq\,\B^n\times\B^n$ as $\to_{\ga}\;=\;\to_{\e}$.
\end{definition}

\begin{figure}[t!]
	\centerline{
  \begin{minipage}{.41\textwidth}
	  \centerline{\scalebox{.85}{\begin{tikzpicture}[>=to,auto]
			\tikzstyle{type} = []
			\tikzstyle{conf} = [rectangle, draw, thick]
			\tikzstyle{pf} = [rectangle, draw, thick, fill=black!10]
			\tikzstyle{lc} = [rectangle, draw, thick, fill=black!70]
			\node[conf](n000) at (0,0) {$\zero\zero\zero$};
			\node[conf](n001) at (2,0) {$\zero\zero\one$};
			\node[conf](n010) at (0,2) {$\zero\one\zero$};
			\node[pf](n011) at (2,2) {$\zero\one\one$};
			\node[pf](n100) at (1,1) {$\one\zero\zero$};
			\node[conf](n101) at (3,1) {$\one\zero\one$};
			\node[conf](n110) at (1,3) {$\one\one\zero$};
			\node[conf](n111) at (3,3) {$\one\one\one$};
			\draw[thick, ->] (n000) edge (n100);
			\draw[thick, ->] (n000) edge[loop left,distance=4mm] (n000);
			\draw[thick, ->] (n000) edge (n001);
			\draw[thick, ->] (n100) edge[loop left,distance=4mm] (n100);
			\draw[thick, ->] (n001) edge[loop right,distance=4mm] (n001);
			\draw[thick, ->] (n001) edge (n011);
			\draw[thick, ->] (n110) edge[loop left,distance=4mm] (n110);
			\draw[thick, ->] (n110) edge (n100);
			\draw[thick, ->] (n010) edge (n110);
			\draw[thick, ->] (n010) edge (n000);
			\draw[thick, ->] (n010) edge (n011);
			\draw[thick, ->] (n011) edge[loop right,distance=4mm] (n011);
			\draw[thick, ->] (n111) edge (n011);
			\draw[thick, ->] (n111) edge (n101);
			\draw[thick, ->] (n111) edge (n110);
			\draw[thick, ->] (n101) edge (n001);
			\draw[thick, ->] (n101) edge[loop right,distance=4mm] (n101);
			\draw[thick, ->] (n101) edge (n100);
		\end{tikzpicture}}}
  \end{minipage}
  \hfill\vrule\hfill
  \begin{minipage}{.53\textwidth}
	  \centerline{\scalebox{.85}{\begin{tikzpicture}[>=to,auto]
\path[use as bounding box] (-0.7,-0.23) rectangle (3.7,3.55);
			\tikzstyle{type} = []
			\tikzstyle{conf} = [rectangle, draw, thick]
			\tikzstyle{pf} = [rectangle, draw, thick, fill=black!10]
			\tikzstyle{lc} = [rectangle, draw, thick, fill=black!70]
			\node[conf](n000) at (0,0) {$\zero\zero\zero$};
			\node[conf](n001) at (2,0) {$\zero\zero\one$};
			\node[conf](n010) at (0,2) {$\zero\one\zero$};
			\node[pf](n011) at (2,2) {$\zero\one\one$};
			\node[pf](n100) at (1,1) {$\one\zero\zero$};
			\node[conf](n101) at (3,1) {$\one\zero\one$};
			\node[conf](n110) at (1,3) {$\one\one\zero$};
			\node[conf](n111) at (3,3) {$\one\one\one$};
			\draw[thick, ->] (n000) edge[loop left,distance=4mm] (n000);
			\draw[thick, ->] (n100) edge[loop left,distance=4mm] (n100);
			\draw[thick, ->] (n001) edge[loop right,distance=4mm] (n001);
			\draw[thick, ->] (n101) edge[loop right,distance=4mm] (n101);
			\draw[thick, ->] (n110) edge[loop left,distance=4mm] (n110);
			\draw[thick, ->] (n011) edge[loop right,distance=4mm] (n011);

		\draw [thick, ->] (n000) edge (n100);
		\draw [thick, ->] (n000) edge (n001);
        \draw [thick, ->] (n000) edge[bend left=5] (n101);
		\draw [thick, ->] (n001) edge (n011);
		\draw [thick, ->] (n101) edge (n001);
		\draw [thick, ->] (n101) edge (n100);
		\draw [thick, ->, bend left=10] (n101) edge (n000);
		\draw [thick, ->] (n110) edge (n100);
		\draw [thick, ->] (n010) edge (n110);
        \draw [thick, ->] (n010) edge (n000);
		\draw [thick, ->] (n010) edge (n011);
		\draw [thick, ->] (n010) edge (n100);
		\draw [thick, ->, bend left=10] (n010) edge (n111);
		\draw [thick, ->, bend right=40, looseness=1.25] (n010) edge (n001);
        \draw [thick, ->] (n010) edge[bend left=5] (n101);
		\draw [thick, ->, bend right, looseness=1.2] (n111) edge (n100);
        \draw [thick, ->] (n111) edge (n101);
		\draw [thick, ->, bend right=70, looseness=1.75] (n111) edge (n000);
		\draw [thick, ->] (n111) edge (n001);
		\draw [thick, ->] (n111) edge (n110);
		\draw [thick, ->] (n111) edge (n010);
		\draw [thick, ->, bend right] (n111) edge (n011);
		\end{tikzpicture}}}
  \end{minipage}}
	\caption{Fully-asynchronous (left) and asynchronous (right) dynamics of BN $f$ defined in Example~\ref{ex:Gf}.}
	\label{fig:fa_a}
\end{figure}
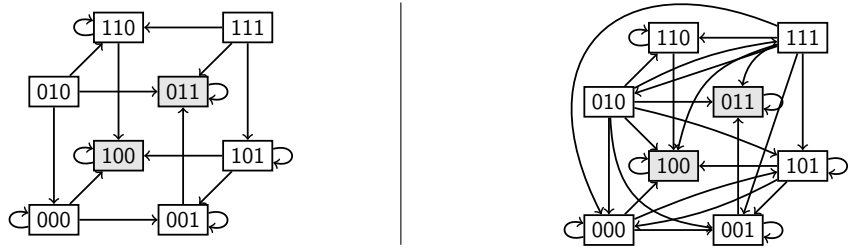

One of the most usual non-deterministic updating modes of BNs is the \emph{fully-asynchronous} mode\footnote{The \emph{fully-asynchronous} mode is usually referred to as \emph{asynchronous} in the system biology modeling community.}, where only one automaton is updated in a transition.
It is largely employed for the analysis of models of biological systems, arguing it enables capturing (some) behaviors caused by different time scale for automata updates.
\begin{definition}
    The \emph{fully-asynchronous} updating mode of a BN $f$ generates the transition relation $\to_{\fa}\ \subseteq\,\B^n\times\B^n$ such that, for all configurations
    $x,y\in\B^n$: $x \to_{\fa} y$ if and only if there exists $i\in\range n$ with $y=\phi_i(x)$.
\end{definition}

\autoref{fig:fa_a} shows the dynamics generated by the fully-asynchronous and asynchronous updating
modes on the BN of \autoref{ex:Gf}.

\section{Non-deterministic updates as set updates}

The updates considered so far are deterministic, and can thus be defined as functions mapping configurations, i.e., of the form $\phi: \B^n\to \B^n$.
As we have seen above, deterministic updates can generate non-deterministic updating modes, by allowing different updates to be applied on a same configuration.

Let us now extend to non-deterministic updates, that we model by functions mapping \emph{sets of configurations}, i.e., of the form $\Phi: 2^{\B^n}\to 2^{\B^n}$.
We define $\Phi$ as a map from sets of configurations to sets of configurations for enabling iterations and compositions of non-deterministic updates.
Nevertheless, we assume that for any $X\subseteq \B^n$, $\Phi(X) = \bigcup_{x\in X} \Phi(\{x\})$:
one can define $\Phi$ only from all singleton configuration set.
This restriction ensures that, for any $X\subseteq\B^n$, each configuration in the image set $y\in\Phi(X)$ can be computed from a singleton set $\{x\}$ for some $x\in\B^n$.
In the following, we call such updates \emph{set updates}.

Starting from a singleton configuration set $\{x\}$, the iteration of set updates delineates the domains of configurations the system can evolve to.
Thus, set updates naturally define transition relations between configurations:
\begin{definition}
    Given a set update function $\Phi$ for BNs of dimension $n$, the generated transition relation is given by $\delta: (2^{\B^n}\to 2^{\B^n}) \to
    2^{\B^n\times\B^n}$ with
    $\delta(\Phi) = \{ (x,y) \mid x\in\B^n, y\in\Phi(\{x\})\}$.
\end{definition}

In contrast with deterministic updates, non-deterministic updating modes can be characterized
directly by set updates.
Indeed, 
non-deterministic updating modes allow ``superposing'' alternative updates to generate different transitions from a single configuration $x$, although each of them is computed with a deterministic update.
For instance, with one update $\phi$ where $\phi(x) = y$ and another update $\phi'$ where $\phi'(x)=y'\neq y$.
Now, let us imagine an updating mode superposing two set updates, $\Phi$ and $\Phi'$ where, for some configurations $x\in\B^n$, $\Phi(\{x\})\setminus\Phi'(\{x\})\neq \emptyset$.
One can then build a single set update $\Phi^*$ such that $\Phi^*(X)=\Phi(X)\cup\Phi'(X)$.
It results that $\delta(\Phi^*)=\delta(\Phi)\cup\delta(\Phi')$, thus the updating mode can be assimilated to $\Phi^*$.

Finally, notice that limit sets of the generated dynamics $\delta(\Phi)$ can be characterized
as the $\subseteq$-smallest sets of configurations $X\subseteq \B^n$ such that $\Phi(X)=X$.

\section{Updating modes selecting elementary transitions}

With deterministic updates as building blocks, we have seen that one can define non-deterministic updating modes by superposing different update functions. 
The resulting transition relation is then the union of the transition relation generated by each individual update (each of them giving a deterministic
dynamics).
Set updates offer an alternative way to formalize the resulting dynamics, by directly defining the set of out-going transitions from a given configuration.
As we will illustrate with the memory updating mode below, this enables a fine-grained selection of the elementary transitions which may then depend on the configuration.

\subsection{Asynchronous and fully-asynchronous updating modes}
\label{sec:set-async}

As a first illustration of set updates and how they can characterize updating modes,
consider the following set update for BNs of dimension $n$:
\[
\Phi_{\e}(X) = \{\phi_{W}(x)\mid x\in X, \emptyset\neq W\subseteq \range n\}\text{.}
\]
This set update generates exactly all the elementary transitions: $\delta(\Phi_{\e}) =\; \to_{\e}$.
Thus, $\Phi_{\e}$ characterizes the asynchronous updating mode.
Similarly, let us now consider the following set update:
\[
    \Phi_{\fa}(X) = \{\phi_i(x)\mid x\in X, i\in\range n\}\text{.}
\]
Remark that $\delta(\Phi_{\fa})=\to_{\fa}$, i.e., $\Phi_{\fa}$ characterizes the fully-asynchronous
updating mode.

\subsection{Memory updating mode}

Until now, all the updating modes that have been discussed depend on deterministic updates that are context free, which leads to deal with memoryless dynamical systems.
In~\cite{J-Goles2020,C-Goles2019} have been introduced another model of BNs, called Memory Boolean networks (MBNs).
The first objective of MBNs is to capture the biologically relevant gene-protein BN model introduced
in~\cite{C-Graudenzi2009}, that builds on the following principles:
\begin{itemize}
\item automata are split in two types: a half models genes, the other half models their associated one-to-one proteins;
\item each protein has its own decay time: the number of time steps during which it remains present in the cell after having been produced by the punctual expression of its associated gene.
\end{itemize}

In their original definition given below, MBNs of dimension $n$ are BNs of dimension $n$ parameterized with a vector $\M\in\N_{>0}^n$, setting the maximal delay (called memory) for the degradation of each automaton. 
Then, an automaton is considered active (Boolean $\one$) whenever its delay to degradation is not $0$. 
Formally, MBN are defined as follows:

\begin{definition}
    A \emph{Memory Boolean network} of dimension $n$ is the couple of a BN $f$ of dimension $n$
    and of a memory vector $\M = (\M_1, \dots, \M_n) \in \N_{>0}^n$.
	The set of its configurations is defined as $X_{(f,\M)} = \{(x,d) \in \B^n \times \N^n \mid \forall i \in \range n,\ d_i \in \{0, \dots, \M_i\},\ x_i = 0\ \iff\ d_i = 0 \text{ and } x_i = 1\ \iff\ d_i \in \{1, \dots, \M_i\}\}$.
	The dynamical system $((f,\M),\p)$ is defined by the transition graph $\D_{((f,\M),\p)}$, with $\p$ the parallel updating mode, made of transitions based on updating function $\phi^\star: X_{(f,\M)} \to X_{(f,\M)}$ depending on the memories such that:
	\begin{equation*}
		\forall (x,d), (y,d') \in X_{(f,\M)},\ (x, d) \trans{}_{((f,\M),\p)} (y,d')\ \iff\ (y, d') = \phi_{\range n}^\star(x, d)\text{,}
	\end{equation*}
	where $\forall i \in \range n,\ \phi_{\range n}^\star(x, d)_i = (y_i, d'_i)$, 
	with:
	\begin{equation*} 
		d'_i = \begin{cases}
			0 & \text{if } f_i(x) = 0 \text{ and } d_i = 0\text{,}\\
			d_i - 1 & \text{if } f_i(x) = 0 \text{ and } d_i \geq 1\text{,}\\
			\M_i & \text{if } f_i(x) = 1\text{,}
		\end{cases}
		\quad \text{ and } \quad 
		y_i = \begin{cases}
			1 & \text{if } d'_i \geq 1\text{,}\\
			f_i(x) & \text{if } d'_i = 0\text{.}
		\end{cases}
	\end{equation*}
\end{definition}
From this initial definition, it is easy to see that the dynamics of a MBN is deterministic and operates on discrete configurations that are not Boolean anymore.
But we will see that MBNs enable to develop a new updating mode, called the memory updating mode, that operates directly on Boolean configurations.

First, let us define $\alpha(x)$ the set of memory configurations corresponding to any binary
configuration $x\in\B^n$,
and conversely, $\beta(d)$ the binary configuration corresponding to a memory configuration
$d\in\N^n$.
Notice that $\forall x\in \B^n$, $\forall d\in\alpha(x)$, $\beta(d)=x$.
\begin{align*}
    \alpha(x)&=\{d\in \mathbb N^n\mid x_i=0 \Leftrightarrow d_i=0, x_i=1\Leftrightarrow d_i\in\range{\M_i}\}\text,
    \\
    \forall i\in\range n\quad \beta(d)_i &= \min \{d_i, 1\}\text{.}
\end{align*}

It appears that $X_{(f,\M)} = \{ (\beta(d),d) \mid d\in\N^n, \forall i\in\range n,
d_i\in\{0,\dots,M_i\} \}$. Thus one can reformulate the
original definition by considering the deterministic parallel update of memory configurations
$d\in\N^n$, and replacing $x$ with $\beta(d)$:
an automaton $i\in\range n$ is set to state $\M_i$ whenever its local function $f_i$ is evaluated to
$\one$ on the corresponding binary configuration $\beta(d)$;
otherwise, its state is decreased by one, unless it is already $0$.
In particular, one can define the  deterministic memory update $\phi^\ast_{\M}:\N^n\to\N^n$ such
that, for each $i\in\range n$,
\[
			\phi^\ast_{\M}(d)_i = \begin{cases}
				0 &\text{if }f_i(\beta(d))=\zero \land d_i=0\text{,}\\
				d_i-1 &\text{if }f_i(\beta(d))=\zero\land d_i\geq 0\text{,}\\
				\M_i &\text{if }f_i(\beta(d))=\one\text{.}
			\end{cases}
\]

Let us now extend the above definitions to sets:
\begin{itemize}
    \item $\forall X\subseteq \B^n,\, A(X)=\bigcup_{x\in X}\alpha(x)$;
    \item $\forall D\subseteq \N^n,\, B(D)=\{\beta(d)\mid d\in D\}$;
    \item $\forall D\subseteq \N^n,\, \Phi^\ast_{\M}(D)=\{\phi^*_{\M}(d)\mid d\in D\}$.
\end{itemize}

The memory set update can then be defined for any set of configurations $X\subseteq \B^n$ by first
generating the set of corresponding memory configurations, then applying the deterministic update on them,
and finally converting them back to binary configurations:
\[
    \Phi_{\M}(X) = B \circ \Phi^\ast_{\M}\circ A (X)\text{.}
\]

With this formulation, one can see that the memory updating mode, being the projection of MBN
configurations on their binary part, lead to non-deterministic dynamics.
Indeed, whenever a configuration gets mapped to several possible memory configurations, and whenever for two of these configurations $d$ and $d'$, there is an automaton $i\in\range n$ where $\phi_{\M}(d)_i=0$ and $\phi_{\M}(d')_i\geq 1$.
This can occur if and only if $\M_i\geq 2$, $x_i=1$, and $f_i(x)=0$.
Thus, the memory updating mode of BNs can equivalently be parameterized by a set of automata $\Mb = \{i\in\range n\mid \M_i\geq 2\}$ and defined as the following set update:
\[
    \Phi_{\Mb}(X) = \{ \phi_W(x)\mid x\in X,
        W\subseteq \range n,
    W\supseteq \{i\in\range n\mid i\notin{\Mb}\vee f_i(x)=\one\}\}\text.
\]
Remark that this definition no longer relies on memory configurations in $\N^n$.
Overall, the memory updating mode of BNs can be understood as a particular set of
elementary transitions: those where automata not in $\Mb$ or automata that can change from state 0
to 1 are always updated, together with any subset of the others (automata in $\Mb$ that can change from
state 1 to 0):
automata in $\Mb$ that are decreasing are updated asynchronously, while the others are updated in
parallel.

\begin{figure}[t!]
	\centerline{
  \begin{minipage}{.25\textwidth}
	  \centerline{\scalebox{.85}{\begin{tikzpicture}[>=to,auto]
\path[use as bounding box] (-0.6,-0.55) rectangle (3.6,3.55);
			\tikzstyle{type} = []
			\tikzstyle{conf} = [rectangle, draw, thick]
			\tikzstyle{pf} = [rectangle, draw, thick, fill=black!10]
			\tikzstyle{lc} = [rectangle, draw, thick, fill=black!70]
			\node[conf](n000) at (0,0) {$\zero\zero\zero$};
			\node[conf](n001) at (2,0) {$\zero\zero\one$};
			\node[conf](n010) at (0,2) {$\zero\one\zero$};
			\node[pf](n011) at (2,2) {$\zero\one\one$};
			\node[pf](n100) at (1,1) {$\one\zero\zero$};
			\node[conf](n101) at (3,1) {$\one\zero\one$};
			\node[conf](n110) at (1,3) {$\one\one\zero$};
			\node[conf](n111) at (3,3) {$\one\one\one$};
            \draw[thick, ->] (n000) edge[bend left=5] (n101);
			\draw[thick, ->] (n101) edge (n100);
		\draw [thick, ->, bend left=10] (n101) edge (n000);
			\draw[thick, ->] (n100) edge[loop left,distance=4mm] (n100);
			\draw[thick, ->] (n110) edge (n100);
			\draw[thick, ->] (n010) edge (n101);
			\draw[thick, ->] (n011) edge[loop right,distance=4mm] (n011);
		\draw [thick, ->, bend right, looseness=1.2] (n111) edge (n100);
		\draw [thick, ->, bend right=70, looseness=1.75] (n111) edge (n000);
			\draw[thick, ->] (n001) edge (n011);
		\end{tikzpicture}}}
  \end{minipage}
  \hfill\vrule\hfill
  \begin{minipage}{.25\textwidth}
	  \centerline{\scalebox{.85}{\begin{tikzpicture}[>=to,auto]
\path[use as bounding box] (-0.6,-0.55) rectangle (3.6,3.55);

			\tikzstyle{type} = []
			\tikzstyle{conf} = [rectangle, draw, thick]
			\tikzstyle{pf} = [rectangle, draw, thick, fill=black!10]
			\tikzstyle{lc} = [rectangle, draw, thick, fill=black!70]
			\node[conf](n000) at (0,0) {$\zero\zero\zero$};
			\node[conf](n001) at (2,0) {$\zero\zero\one$};
			\node[conf](n010) at (0,2) {$\zero\one\zero$};
			\node[pf](n011) at (2,2) {$\zero\one\one$};
			\node[pf](n100) at (1,1) {$\one\zero\zero$};
			\node[conf](n101) at (3,1) {$\one\zero\one$};
			\node[conf](n110) at (1,3) {$\one\one\zero$};
			\node[conf](n111) at (3,3) {$\one\one\one$};

    \begin{scope}[densely dotted]
		\draw [->] (n000) edge (n100);
		\draw [->] (n000) edge (n001);
        \draw [->] (n000) edge[bend left=5] (n101);
		\draw [->] (n001) edge (n011);
		\draw [->] (n101) edge (n001);
		\draw [->] (n101) edge (n100);
		\draw [->, bend left=10] (n101) edge (n000);
		\draw [->] (n110) edge (n100);
		\draw [->] (n010) edge (n110);
        \draw [->] (n010) edge (n000);
		\draw [->] (n010) edge (n011);
		\draw [->] (n010) edge (n100);
		\draw [->, bend left=10] (n010) edge (n111);
		\draw [->, bend right, looseness=1.25] (n010) edge (n001);
        \draw [->] (n010) edge[bend left=5] (n101);
		\draw [->, bend right, looseness=1.2] (n111) edge (n100);
        \draw [->] (n111) edge (n101);
		\draw [->, bend right=70, looseness=1.75] (n111) edge (n000);
		\draw [->] (n111) edge (n001);
		\draw [->] (n111) edge (n110);
		\draw [->] (n111) edge (n010);
		\draw [->, bend right] (n111) edge (n011);
    \end{scope}

		\draw [thick, ->, bend right=70, looseness=1.75] (n000) edge (n111);
        \draw [thick, ->] (n101) edge[bend left=5] (n010);

		\end{tikzpicture}}}
  \end{minipage}
  \hfill\vrule\hfill
  \begin{minipage}{.25\textwidth}
	  \centerline{\scalebox{.85}{\begin{tikzpicture}[>=to,auto]

\path[use as bounding box] (-0.6,-0.55) rectangle (3.6,3.55);

			\tikzstyle{type} = []
			\tikzstyle{conf} = [rectangle, draw, thick]
			\tikzstyle{pf} = [rectangle, draw, thick, fill=black!10]
			\tikzstyle{lc} = [rectangle, draw, thick, fill=black!70]
			\node[conf](n000) at (0,0) {$\zero\zero\zero$};
			\node[conf](n001) at (2,0) {$\zero\zero\one$};
			\node[conf](n010) at (0,2) {$\zero\one\zero$};
			\node[pf](n011) at (2,2) {$\zero\one\one$};
			\node[pf](n100) at (1,1) {$\one\zero\zero$};
			\node[conf](n101) at (3,1) {$\one\zero\one$};
			\node[conf](n110) at (1,3) {$\one\one\zero$};
			\node[conf](n111) at (3,3) {$\one\one\one$};

    \begin{scope}[densely dotted]
		\draw [->] (n000) edge (n100);
		\draw [->] (n000) edge (n001);
        \draw [->] (n000) edge[bend left=5] (n101);
		\draw [->] (n001) edge (n011);
		\draw [->] (n101) edge (n001);
		\draw [->] (n101) edge (n100);
		\draw [->, bend left=10] (n101) edge (n000);
		\draw [->] (n110) edge (n100);
		\draw [->] (n010) edge (n110);
        \draw [->] (n010) edge (n000);
		\draw [->] (n010) edge (n011);
		\draw [->] (n010) edge (n100);
		\draw [->, bend left=10] (n010) edge (n111);
		\draw [->, bend right, looseness=1.25] (n010) edge (n001);
        \draw [->] (n010) edge[bend left=5] (n101);
		\draw [->, bend right, looseness=1.2] (n111) edge (n100);
        \draw [->] (n111) edge (n101);
		\draw [->, bend right=70, looseness=1.75] (n111) edge (n000);
		\draw [->] (n111) edge (n001);
		\draw [->] (n111) edge (n110);
		\draw [->] (n111) edge (n010);
		\draw [->, bend right] (n111) edge (n011);
    \end{scope}

		\draw [thick, ->] (n000) edge (n110);
		\draw [thick, ->, bend right=70, looseness=1.75] (n000) edge (n111);
        \draw [thick, ->] (n000) edge[bend left=15] (n010);
		\draw [thick, ->, bend right=20, looseness=2.00] (n000) edge (n011);
		\draw [thick, ->, bend right] (n101) edge (n110);
        \draw [thick, ->] (n101) edge[bend left=5] (n010);
        \draw [thick, ->] (n101) edge[bend right=15] (n111);
		\draw [thick, ->] (n101) edge (n011);
		\end{tikzpicture}}}
  \end{minipage}}
  \caption{Memory dynamics with $\Mb=\{1\}$ (left),
  interval dynamics (center), and 
  MP dynamics (right) of the BN of \autoref{ex:Gf}.
  In these two latter, the elementary transitions are dotted
  and loops are omitted.}
	\label{fig:und}
\end{figure}
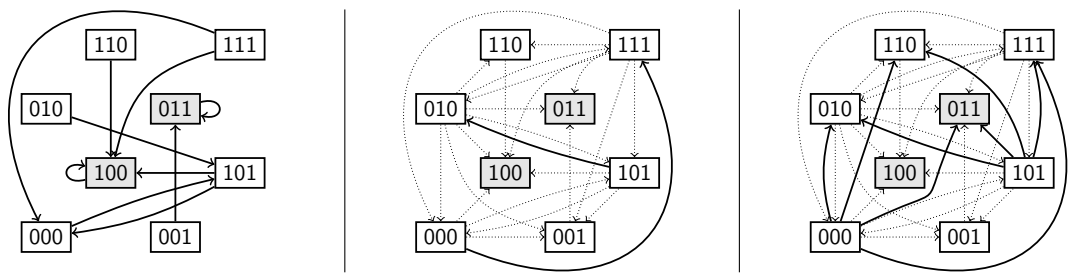

\autoref{fig:und}(left) gives the dynamics generated by the interval updating mode on the BN of
\autoref{ex:Gf} with $\Mb=\{1\}$.

\section{Updating modes going beyond (non-)elementary transitions}

BNs are widely used to model dynamics of biological systems, notably implying gene regulation.
Gene regulation is a dynamical biological process that involves numerous mechanisms and entities
among which some of them, like RNAs and proteins, have specific influences that depend on their
concentration. In other terms, the regulation process in its whole admits significant quantitative
parts.
The question arises then of how faithful are Boolean dynamics with respect to the quantitative dynamics.
It has been recently underlined in~\cite{MPBNs} that the elementary and non-elementary transitions
of BNs are not complete enough to capture particular quantitative trajectories. With a fixed logic,
and starting from similar
configurations, the quantitative system shows that an automaton can eventually get activated, whereas the
asynchronous dynamics of the BN shows it is impossible.

In this section, we address the set update reformulation of two recently introduced dynamics of BNs
which generate transitions that are neither elementary nor non-elementary: they
result in set updates $\Phi$ where, for some BNs of dimension $n$ and for some
configurations $x\in\B^n$,
there is $k\in\N$ such that there exists $y\in\Phi^k(\{x\})$ whereas $x\not\to_{\e}^* y$.

\subsection{Interval updating mode}

From the concurrency theory, it is known that the execution of 1-bounded contextual
Petri nets corresponding to the asynchronous updating mode (known as \emph{steps} semantics) may miss some
transitions that can be triggered when considering certain delay of state
change~\cite{interval-CPNs,bn-concurrency}.
In~\cite{beyond-general} is proposed a translation of the \emph{interval} semantics of Petri nets to
BNs, showing it can predict transitions that are neither elementary nor non-elementary transitions: configurations that are not reachable with the asynchronous mode from a fixed initial
configuration $x$ become reachable with \emph{Interval Boolean networks}.

The main principle of the interval dynamics is to decompose the change of state of an automaton,
allowing interleaving the update of other automata. For instance, let us assume that automaton $i$
can change from state $0$ to $1$. In the interval dynamics, we register that $i$ will eventually
change to state $1$, and then allow the update of other automata, still considering that $i$ is in
state $0$.
In~\cite{beyond-general}, this is applied to any BN of dimension $n$ by an encoding
with the fully-asynchronous dynamics of a BN of dimension $2n$:
each automaton is split in a \emph{read} and \emph{write} automaton,
where the write automaton register the next state of the original automaton,
and the read automaton keeps its current state and will eventually copy the value of the write
automaton.
The dynamics of the original BN is then obtained by projecting the configurations on the read automata.

We provide below an equivalent formulation as a composition of set updates.
Essentially, whenever an automaton $i$ can change its state, we hold it and compute the possible
state changes of the automata different than $i$.
Thus, during this evaluation, we have a growing number of held automata, that we denote by
$L$, waiting for their state to be updated.
The set update $\PhiI L$ extends a given set of configurations with the possible state change of
automata not in $L$.
The function $\PhiD L i(x)$ first computes all possible state changes by iterating $\PhiI {L\cup
\{i\}}$ until a fixed point, and then apply the state change for the automaton $i$ on all the resulting
configurations.
\begin{definition}
    The \emph{Interval set update} $\Phi_{\interval}$ of a BN of dimension $n$ is given by
    $\Phi_{\interval} = \PhiI{\emptyset}$ where
\begin{eqnarray*}
    \PhiI L (X) &=& X\cup 
    \{ y \in \PhiD L i (x) \mid x\in X,
    i\in\range n, i\notin L, f_i(x)\neq x_i\}\text,
    \\
    \PhiD L i (x) &=& 
        \{ \flip y i\mid y\in\PhiI{L\cup\{i\}}^\omega(\{x\})\}\text.
\end{eqnarray*}
\end{definition}
The interval updating mode preserves the fixed points of $f$: for any configuration $x\in\B^n$,
$\Phi_{\interval}(\{x\}) = \{x\}$ if and only if $f(x)=x$.
Moreover, one can prove that it includes all the elementary transitions:
for any configuration $x\in\B^n$, $\Phi_{\e}(\{x\})\subseteq\Phi_{\interval}(\{x\})$.

\begin{example}
    \autoref{fig:und} shows the transitions generated by the interval updating mode on the
    BN of \autoref{ex:Gf}.
    Notice that there is a path from $000$ to $111$, which does not exist in the asynchronous
    dynamics (\autoref{fig:fa_a}).
    Indeed, let us partially compute
    $\PhiI \emptyset (\{000\}) = \{000\} \cup \PhiD \emptyset 1 (000) \cup \PhiD\emptyset 3 (000)$.

Let us focus on the interval update of automaton $1$ with $\PhiD\emptyset 1(000)$, which requires
computing all the iterations of $\PhiI{\{1\}} (\{000\})$.
The first iteration gives
    $\PhiI{\{1\}} (\{000\}) = \{000\}\cup\PhiD{\{1\}} 3(000) = \{000,001\}$,
then
    $\PhiI{\{1\}}^2 (\{000\}) = \{000,001\} \cup \PhiD{\{1\}} 2(001) = \{000,001,011\}
    =\PhiI{\{1\}}^\omega(\{000\})$.

    Finally, we get $\PhiD\emptyset 1(000) = \{100,101,111\}$. Thus,
    $111\in\Phi_{\interval}(\{000\})$.
\end{example}

\subsection{Most Permissive updating mode}
\label{sec:mp}

Most Permissive Boolean networks (MPBNs) have been designed to capture all automata updates that
could occur in any quantitative refinement of the BN.
We will come back more formally to this notion later in this section.

The main feature of MPBNs is to abstract all the possible interaction thresholds between automata.
Consider the case whenever the state of an automaton $i$ is used to compute the state of two
distinct automata $j$ and $k$, and assume that $i$ is increasing from $0$:
during its increase, there are times when $i$ may be high enough for trigger a state
change of $j$ but not (yet) high enough for $k$.
This can be illustrated on a concrete biological example, the so-called \emph{incoherent
feed-forward loop of type 3}~\cite{Mangan03}:
a BN $f$ of dimension 3 with
\[f_1(x)=1 \qquad f_2(x)=x_1 \qquad f_3(x)=\neg x_1\wedge x_2\text.\]
\label{eq:i3ffl}
Starting from the configuration $000$, the asynchronous updating mode predicts only the following
non-reflexive transitions:
$000\to_{\ga} 100 \to_{\ga} 110$.
Notice that in this case, the interval updating mode results in the same transitions.
However, it has been observed experimentally~\cite{Schaerli2014} and in quantitative
models~\cite{Ishihara2005,Rodrigo2011} that depending on reaction kinetics, one can actually
activate transiently the automata $3$.
Essentially, the idea is that during the increase of the state of automaton $1$, there is period
of time where $1$ is high enough so $2$ can consider it active ($x_1$ true) but
$3$ still considers it inactive ($x_1$ false). Then, the state of automaton $2$ can increase, and so
do the state of automaton 3.
This activation of 3 cannot be predicted with BN updating modes defined so far, whereas the logic
encoded by $f$ is correct.

Without introducing any parameter,
MPBNs capture these additional dynamics by accounting for all possible thresholds ordering,
for all updates that can happen between a switch of a Boolean state. In some sense, the MP updating
mode abstracts both the quantitative domain of automata and the duration of state changes.
Their original definition~\cite{MPBNs} is based on the introduction pseudo \emph{dynamic} states,
namely increasing and decreasing.
An automaton can change from 0 to increasing whenever it can interpret the state of the other
automata so that its local function is satisfied. Once in increasing state, it can change to the
state 1 without any condition, or to the decreasing state whenever it can interpret the state of
other automata so that its local function is not satisfied.
Whenever an automaton is in a dynamic state, the automata can freely interpret its state as either 0
or 1.
Remark that the possible interpretations of the MP configurations always result in a hypercube (a
set of automata fixed to a Boolean value, and the others free).

Here, we show that the MP dynamics can be expressed in a more standard way by the means of
composition of set updates.
A first stage consists in \emph{widening} all the elementary set updates to compute all the
possible interpretations of automata changing of state.
The widening is defined using the function $\nabla:2^{\B^n}\to2^{\B^n}$ which computes the vertices of
the smallest hypercube containing the given set of configurations.
For instance, $\nabla(\{01,10\})=\{00,01,10,11\}$.
Given a set of automata $W$,
the widening set update $\Phi_{W, \nabla}:2^{\B^n}\to2^{\B^n}$ applies this operator on the results
of the elementary set update, or equivalently with the fully-asynchronous set update, on the
automata of $W$ (\autoref{sec:set-async}).
This widening is re-iterated until a fixed point is reached.
Then, a \emph{narrowing} $\Lambda_W: 2^{\B^n}\to 2^{\B^n}$ filters the computed configurations $X$ to retain only those where the states
of automata in $W$ can be computed with $f$ from $X$.
\begin{definition}
    The \emph{Most Permissive set update} $\Phi_{\MP}$ of a BN of dimension $n$ is given by
\[\Phi_{\MP}(X) = \bigcup_{W\subseteq\range n} \Lambda_W\circ \Phi^\omega_{W,\nabla}(X)\]
where, for any $X\subseteq\B^n$ and any $W\subseteq\range n$:
\begin{eqnarray}
    \nabla(X) &=& \{ x\in\B^n\mid \forall i\in\range n, \exists y\in X: x_i=y_i\}\text,
    \\
    \Phi_{W,\nabla}(X) &=& \nabla(X \cup \{\phi_i(x)\mid x\in X, i\in W\})\text,
    \\
    \Lambda_W(X) &=& \{x\in X\mid \forall i\in W, \exists y\in X: x_i=f_i(y)\}\text.
\end{eqnarray}
\end{definition}

\begin{example}
\autoref{fig:und} shows the dynamics generated by the MP updating mode on the BN of \autoref{ex:Gf}.
With the BN $f$  of the incoherent feed-forward loop introduced at the beginning \autoref{sec:mp}
page~\pageref{eq:i3ffl}, we obtain:
\begin{eqnarray*}
    \Phi_{\{1,2,3\},\nabla}(\{000\}) &=& \nabla(\{000,100\}) =\{000,100\}\text,\\
    \Phi_{\{1,2,3\},\nabla}^2(\{000\}) &=& \nabla(\{000,100\}\cup\{110\}) =
    \{000,100,010,110\}\text,\\
    \Phi_{\{1,2,3\},\nabla}^3(\{000\}) &=& \nabla(\{000,1000,010,110\}\cup\{011\}) = \B^n\text,\\
    \Lambda_{\{1,2,3\}}(\B^n) &=& \{100,101,110,111\}\text.
\end{eqnarray*}
Thus, $111\in\Phi_{\MP}(\{000\})$, whereas $000\not\to_{\e}^* 111$ and $111\notin\Phi_{\interval}(\{000\})$.
\end{example}

Let us now list some basic properties of the MP updating mode:
\begin{enumerate}
    \item MP preserves the fixed points of $f$: for any configuration $x\in \B^n$, $f(x)=x$ if and only if $\Phi_{\MP}(x) = \{x\}$.
    \item MP subsumes elementary transitions: $\to_{\e}\,\subseteq\, \delta(\Phi_{\MP)}$.
    \item MP transition relation is transitive and reflexive: $\Phi_{\MP} = \Phi_{\MP}^2$.
    \item (by 2 and 3) MP transition relation subsumes non-elementary transitions:
        $\to^*_{\e}\,\subseteq\, \delta(\Phi_{\MP})$.
    \item (by 4 and the example) there exist BNs $f$ such that the MP transition relation is strictly larger than non-elementary transitions, i.e.,
        there exist $x,y\in\B^n$ such that $y\in\Phi_{\MP}(\{x\})$ but
        $x\not\to_{\e}^* y$.
\end{enumerate}

In~\cite{MPBNs}, it has been demonstrated that MP dynamics of a BN $f$ forms a correct abstraction of the
dynamics of any quantitative model being a \emph{refinement} of $f$.
A quantitative model $F$ can be defined as a function mapping discrete or continuous
configurations to the derivative of the state of automata. Then, $F$ is a refinement of $f$ if and
only if the derivative of automaton $i$ is strictly positive (resp.\ negative) in a given
quantitative configuration $z$ only if there is a binarization $\tilde z$ of $z$ so that $f_i(\tilde
z)=1$ (resp. $0$).
It has also been proven to be minimal for the abstraction of asynchronous discrete models.
Moreover, the complexity for deciding the existence of a path between two configurations as well as
deciding whether a configuration belongs to a limit set is respectively in P$^{\text{NP}}$ and in
coNP$^{\text{coNP}}$ in general and in P and in coNP for locally monotonic BNs (each local
function is monotonic with respect to a specific component-wise ordering of configurations),
in contrast with the other updating modes where these problems are PSPACE-complete.

\section{Discussion}

By extending to non-deterministic updates modeled as set updates, we can
reformulate in a unified manner a range of BN dynamics introduced in the literature with ad-hoc
definitions, and for which the usual deterministic updates seem not expressive enough.
These reformulations bring a better understanding and comparison of dynamics as more classical BN
updating modes.
Moreover, they allow envisioning new families of updating modes as variations of the one presented
here.
For instance, the given MP set update allows to readily define restrictions of it:
similarly to the block-sequential updating mode, one could parameterize the MP set update to only
consider particular sequences of sets of automata to update.
One could also consider different narrowing operators and different manners to compose them with the widening,
with the goal of reducing the set of generated transitions.

On the one hand, these set updates foster the definitions of totally new kinds of updating modes.
On the other hand, they raise the question of a potential upper limit on which transitions could be considered as valid, or at least reasonable.

\paragraph*{On \emph{reasonable} set updates}

Of course, from a purely theoretical standpoint, any set update which is mathematically correct is reasonable but, if we consider set updates in a context of modeling, some constraints need to be taken in account. 
This second standpoint is the one on which is based the following discussion.
Indeed, as evoked in the introduction of this paper, BNs are a classical mathematical model in systems biology. 
They are notably widely used to model genetic regulation networks, in which their use rests for instance on the fact that their limit sets model real observable ``structures'' such as differentiated cellular types (fixed points), or specific biological paces (limit cycles).
In this sense, a basic criterion would be that an updating mode for a BN $f$ is admissible only if the fixed points of $f$ are fixed points of the generated dynamics as well.
This criterion would allow capturing the fundamental property of fixed point stability of dynamical system theory.
For instance, let us consider the set update $\Phi_\top(X) = \B^n$:
clearly, the set of fixed points of the generated dynamical system is always empty, and thus do not include those of $f$ whenever $f$ has at least one fixed point.
Therefore, such a set update does not appear satisfying.

Now, let us discuss about set updates which would give sets larger than MP for some singleton
configuration set $\{x\}$.
First, what about defining a widening operator larger than $\nabla$?
For any set of automata $W$ and for any configuration $x$, remark that
$\Phi_{W,\nabla}^\omega(\{x\})=Y$ is the smallest hypercube containing $x$ verifying for each
automata $i\in W$ that
for any  configuration $y\in Y$, if $f_i(y)\neq x_i$, then there exists a configuration $z\in Y$
with $z_i\neq x_i$.
Thus, an automaton in $W$ is either fixed to its state in $x$, or it has been computed with its local function from at least one configuration from a smaller hypercube.
Therefore, a widening operator $\nabla'$ verifying for some $X\subseteq \B^n$, $\nabla'(X)\supsetneq
\nabla(X)$ implies that the state of at least one configuration is not computed using $f$ on $X$.
Now, what about a less stringent narrowing operator.
Let us consider a configuration $y\in\Phi_{W,\nabla}(\{x\})=Y$ for some set of automata $W$, but
$y\notin\Lambda_W(Y)$.
This implies that there exists an automaton $i\in W$ such that $\forall z\in Y$, $y_i\neq f_i(z)$, i.e., $y_i$ cannot be computed by $f_i$ from $X$.
Overall, a set update giving configuration sets strictly larger than the MP update implies that for
some configurations, the state of at least one automaton is not computed using its local function.

\paragraph*{Simulations by deterministic updates}

A perspective of the work presented in this paper focuses on simulations of BNs evolving with
non-deterministic updates by BNs evolving with deterministic updates.
A first natural way is by following a classical determinization of the dynamics.
Indeed, one can encode any set of configurations in $\B^n$ as one configuration in $\B^{2^n}$.
Let us consider such an encoding $c:2^{\B^n}\to \B^{2^n}$ where,
for all $x\in\B^n$, $c(X)_{x} = 1$ if $x\in X$, otherwise $c(X)_x = 0$
(we slightly abuse notations here, by specifying a vector index by its binary representation).
Now, it is clear that for any set update $\Phi:2^{\B^n}\to 2^{\B^n}$ of a BN $f$ of dimension $n$,
one can define a BN $g$ such that
for all sets of configurations $X\subseteq\B^n$, $g(c(X)) = c(\Phi(X))$.
This encoding is complete in the sense that any transition generated by $\Phi$ is simulated in $(g,\p)$.
But these simulations are nothing else but a brute-force encoding in which we get rid of the transition relation by increasing exponentially the state space.
Moreover, with this deterministic encoding, the structure of the transition relation of $(f,\mu=\Phi)$ is lost, which make much more difficult characterizing dynamical features of $(f,\mu)$ such as its limit sets for instance.

Actually, a fundamental matter here lies in the concept of simulation at stake here:
we are interested in intrinsic simulations which go far beyond the classical concepts of encoding or simulation. 
Indeed, intrinsic simulations aim at conserving dynamical structures in addition to operated computations. 
So, one of the first question to answer would consist in defining formally different kinds of intrinsic simulations.
Nevertheless, firstly, consider the following intrinsic simulation: a dynamical system $(f,\mu)$ simulates another $(g,\mu')$ if $\D_{(g,\mu')}$ is a subgraph of $\D_{(f,\mu)}$.
With this rather simple definition, it is direct to state that, with $\ga$ and $\Mb$ the asynchronous and memory updating modes respectively, for any BN $f$, $(f,\ga)$ simulates $(f,\Mb)$. 
Some natural questions related to BNs updated with memory are the following: 
\begin{itemize}
\item Are there BNs whose dynamics obtained according to $\Mb$ remains deterministic, whatever $\Mb$?
\item If so, what are their properties and what are the equivalent deterministic updating modes?
\end{itemize}
To go further, consider the $\MP$ updating mode.
It is direct that $(f,\mu)$ does not simulate $(f,\MP)$, except for very particular $f$.
Let us now consider a more general intrinsic simulation: a dynamical system $(f,\mu)$ simulates another $(g,\mu')$ if $\D_{(g,\mu')}$ is a graph obtained from $\D_{(f,\mu)}$ thanks to edge deletions, and vertex shortcuts.
A lot of promising questions arise from this, in particular related to $\Mb$ and $\MP$ updating modes, among which for instance:
\begin{itemize}
\item Let $\per$ be a deterministic periodic updating mode. How can $(f,\Mb)$ be simulated by $(g,\per)$? The answer is known for $\per = \p$~\cite{J-Goles2020}, but it seems pertinent to find a generalization to deterministic periodic updating modes, and even more general deterministic updating modes. 
\item Intuitively, any $(f,\MP)$ might be simulated by $(g,\ga)$, where $f$ and $g$ are BNs and the dimension of $g$ is greater than that of $f$. But how many automata need to be added to $g$ depending on the dimension of $f$?
\end{itemize}
All answers, even partial or negative, will bring a better understanding of updating modes and BNs, which would lead to pertinent further development in both BN theory and their application in systems biology.


\bibliography{automata21.bib}

\end{document}